\documentclass [12pt,oneside] {article}
\usepackage{latexsym,amsfonts,amssymb,amscd,amsmath,amsthm,mathrsfs,stmaryrd} 
\usepackage{fullpage}
\usepackage{lscape}
\usepackage{stmaryrd}
\usepackage{graphicx}

\newcommand {\mF} {\mathcal F}
\newcommand {\mG} {\mathcal G}

\newcommand {\mL} {\mathcal L}
\newcommand {\mM} {\mathcal M}
\newcommand {\mN} {\mathcal N}
\newcommand {\mO} {\mathcal O}
\newcommand {\mR} {\mathcal R}

\newcommand {\CC} {\mathbb C}

\newcommand {\NN} {\mathbb N}
\newcommand {\NNo} {\mathbb N_{\geq 1}}
\newcommand {\RR} {{\mathbb R}}
\newcommand {\QQ} {{\mathbb Q}}
\newcommand {\RZ} {{\RR\mmod\ZZ}}
\renewcommand {\SS} {{\mathbb S}} 

\newcommand {\ZZ} {\mathbb Z}

\newcommand{\map}[3] {#1\negmedspace:#2\to #3}
\newcommand{\longmap}[3] {#1\negmedspace:#2\longrightarrow #3}
\newcommand{\on}[1] {{\operatorname{#1}}}

\newcommand {\at}[1] {\arrowvert_{#1}}

\newcommand {\cen} {{\on{Center}}}

\newcommand {\id} {\on{Id}}
\newcommand {\inv} {^{-1}}
\newcommand {\mmod} {/\!\!/}
\newcommand {\ob} {\on{ob}}
\newcommand {\ol}[1] {\overline{#1}}
\newcommand {\ps}[1] {{[\![#1]\!]}}
\newcommand {\pt} {\on{pt}}

\newcommand {\sub} {\subseteq}
\newcommand {\tensor} {\otimes}
\newcommand {\ul}[1]{\underline{#1}}

\renewcommand {\(} {\left(}
\renewcommand {\)} {\right)}

\newtheorem {Thm} {Theorem}[section]
\theoremstyle{plain}
\newtheorem {Cor}[Thm] {Corallary}
\newtheorem {Lem}[Thm] {Lemma}
\newtheorem {Prop}[Thm] {Proposition}
\theoremstyle{definition}

\newtheorem {Def}[Thm] {Definition}
\newtheorem {Exa}[Thm] {Example}

\usepackage [all]{xy}
\SelectTips{cm}{12}

\newcommand {\gpd} {\mathcal G\medspace\!\!pd}
\newcommand {\fun} {\mathcal F\!un}

\begin{document}

\title{Power operations in orbifold Tate $K$-theory}
\author{Nora Ganter\thanks{The author is supported by a Discovery
    Grant from the Australian Research Council and currently holds an
    Australian Research Fellowship.
  This paper extends work that the author
  pursued in 2006, during a stay at MSRI and the subsequent time at University
  of Illinois. At that time, she was
  supported by NSF grant DMS-0504539.}\\ The University of Melbourne}
\date {\today}

\maketitle
\begin{abstract}
We formulate the axioms of an orbifold theory with power
operations. We define orbifold Tate
$K$-theory, by adjusting Devoto's definition of
the equivariant theory, and proceed to construct its
power operations. We calculate the resulting symmetric powers,
exterior powers and Hecke operators and put our work into context with
orbifold loop spaces, level structures on the Tate curve and
generalized Moonshine.
\end{abstract}
\section{Introduction}
Tate $K$-theory is a form of elliptic cohomology. The corresponding
generalized elliptic curve is the Tate curve, 
and the corresponding
elliptic genus is the Witten genus \cite{Ando:Hopkins:Strickland:eswgtc}. 

The most conceptual definition of the Tate curve is probably the
characterization 
in \cite[Thm VII.2.1]{Deligne:Rapoport}:
if $\mM_{ell}$ is the algebraic stack of (nice enough) generalized elliptic
curves\footnote{We write $\mM_{ell}$ for the stack denoted $\mM_1$ in
  \cite[III.2.6]{Deligne:Rapoport}.} 
then the
completion of $\mM_{ell}$ at infinity,
\begin{eqnarray*}
  spf\medspace\ZZ\ps q & \longrightarrow & \mM_{ell},
\end{eqnarray*}
classifies a formal projective curve. This possesses an 
algebraization, living over $\on{spec}\medspace\!\ZZ\ps q$, which we will call $Tate(q)$.
A more pedestrian definition of $Tate(q)$ can be found in \cite{Ando:Hopkins:Strickland:eswgtc}.

Out of all known elliptic cohomology theories, $K_{Tate}$ might be the
one where the conjectural relationship to string theory is best
understood, and Witten's original definition of his genus involves loop
spaces.

The definition of $G$-equivariant Tate $K$-theory for finite $G$ goes
back to Devoto \cite{Devoto} and is modeled on the loop space of a
global quotient orbifold. Our definition of orbifold Tate $K$-theory
is a mild modification of Devoto's. Devoto's work also already makes
the connection to level structures on elliptic curves.

Guided by the ideas in \cite{Dijkgraaf:Moore:Verlinde:Verlinde}, we
define power operations on $K_{Tate}$. This project grew out of the
desire for a more explicit link between the findings in
\cite{Ganter:thesis} and the string theoretic picture in
\cite{Dijkgraaf:Moore:Verlinde:Verlinde} and was inspired by the
similarities between \cite{Lupercio:Uribe:Xicotencatl},
\cite{Tamanoi} and \cite{Dijkgraaf:Moore:Verlinde:Verlinde}.

The story in \cite{Ganter:thesis} is one of convergence: there, we saw
how the purely homotopy theoretic formalism of power operations in
elliptic cohomology (the theory there is $E_2$) results in a formula
much resembling that of Dijkgraaf, Moore, Verlinde and Verlinde.

Here, we take a more traditional approach to mathematical physics,
modeling the mathematical definitions on the
calculations of physicists. The resulting power operations are closely
linked to level 
structures and isogenies on the Tate curve, suggesting that this is
the `right' definition of elliptic power operations on
$K_{Tate}$. To make precise exactly what is meant by `right' in this
context, one has to revisit the work of Matthew Ando, his power
operations on $K_{Tate}$,
and the definition of {\em elliptic power operations} in
\cite{Ando:Hopkins:Strickland:sigma}. 
All this was pioneered in the work of Andrew Baker, who was the first to
define Hecke operators on elliptic cohomology  
\cite{Baker1}, \cite{Baker2} and \cite{Baker3}.

When restricted to the coefficients of $K_{Tate}$, the total symmetric
power satisfies the identity
\begin{eqnarray}\label{eq:generating}
  S^{Tate}_t(x) & = & \exp\(\sum_{m\geq 1}T_m(x)\),
\end{eqnarray}
where the $T_m$ are the equivariant Hecke operators of generalized
Moonshine defined in \cite{Ganter:moonshine}, \cite{Morava} and prominent in
Carnahan's proof of the generalized Moonshine Conjectures
\cite{Carnahan}. The identity \eqref{eq:generating} gives rise to a
new formulation of the replicability condition, although the latter
remains mysterious. On the other hand, the total symmetric power of a
vector bundle is the
stable exponential characteristic class
\begin{eqnarray}\label{eq:Witten_class}
  S^{Tate}_t(V) & = & \bigotimes_{k\geq 1}S_{t^k}(V),
\end{eqnarray}
used to define the Witten genus.
Here $S_t$ is the total symmetric power in $K$-theory. 
The Witten genus is, of course, an important player in Moonshine,
because of the Hirzebruch conjecture \cite{Hirzebruch:Berger:Jung}
\cite{Hopkins:Mahowald}. 

In the context of replicability, the variable $t$ in \eqref{eq:generating} lives a
peculiar and poorly understood double life as formal variable on one hand and as modular
variable $e^{2\pi i\tau}$ on the other. In the context of
\eqref{eq:Witten_class}, the same double life occurs, and a compelling
explanation for this phenomenon, involving the boundary of a Krichever
style\footnote{Cutting along the circle $[0,1]$ inside
  $\CC/\langle\tau,1\rangle$, instead of Krichever's small circle as
  in \cite[8.11]{Pressley:Segal}.}
construction of the relevant moduli space goes back to unpublished
work of Looijenga \cite{Baranovsky:Ginzburg}.

\subsection{Acknowledgments} It is a pleasure to thank Nitu Kitchloo,
Matthew Ando, Charles Rezk, Alex Ghitza, David Gepner, Craig
Westerland, Christian Haesemeyer, Ralph Cohen, Ernesto Lupercio,
Bernardo Uribe, Jorge Devoto and John McKay for many helpful
conversations on the topic. I would also like to thank the referee for
very useful suggestions on an earlier account of the material and for
pointing out the connection to level structures.
\section{Orbifold Tate $K$-theory}
We will work in the 2-category $\gpd$ whose objects are the
(small) topological groupoids and with 
\begin{eqnarray*}
  1\on{Hom}(X,Y) & = & \fun(X,Y)
\end{eqnarray*}
the groupoid of continuous functors from $X$ to $Y$ and continuous
natural isomorphisms between them.
We do not emphasize the 2-category point of view. For all practical
purposes, we could as well work in the category
\begin{eqnarray*}
  \on{Gpd} & = & \gpd /\cong,
\end{eqnarray*}
obtained from $\gpd$ by identifying naturally isomorphic 1-morphisms, or
even in Moerdijk's orbifold category $Orb$ (see \cite{Moerdijk}).
The {\em center} of a groupoid $X$ is the group 
\begin{equation*}
  \cen(X) : = 2\on{Hom}\(\id_X,\id_X\) = \mN\!at\(\id_X,\id_X\),    
\end{equation*}
of natural transformations from $\id_x$ to $\id_x$.  
We will also need the 2-category $\gpd^{cen}$ whose objects are pairs
$(X,\xi)$ with $\xi$ a center element of $X$ and with 
\begin{eqnarray*}
  1\on{Hom}((X,\xi),(Y,\upsilon)) & \subset & \fun(X,Y)
\end{eqnarray*}
the full subcategory of functors $f$ satisfying
\begin{eqnarray*}
  f\xi & = & \upsilon f.
\end{eqnarray*}
For simplicity of exposition, these center
elements are assumed to be of
finite order. This is automatically the case if $X$ is an orbifold
groupoid in the sense of Moerdijk.
If $G$ is a finite group and $M$ is a $G$-space we write $M\mmod G$
for the corresponding translation groupoid.
\begin{Exa}
  The group
  \begin{eqnarray*}
    \cen(\pt\mmod G) & = & \cen (G)
  \end{eqnarray*}
  is the center of the group $G$.
\end{Exa}
\begin{Exa}
  The inertia groupoid
  \begin{eqnarray*}
    \Lambda X & = & \mF\!un(\pt\mmod \ZZ,X)
  \end{eqnarray*}
  of $X$ admits a group homomorphism
  \begin{eqnarray*}
    \ZZ & \longrightarrow & \cen(\Lambda X)\\
     k  & \longmapsto     & \(\xi^k\negmedspace : f\mapsto f(k)\).
  \end{eqnarray*}
  Viewing objects of $\Lambda X$ as pairs $(x,g)$ with $x\in ob(X)$
  and $g\in aut(x)$, we have
  \begin{eqnarray*}
    \xi^k_{(x,g)} & = & g^k
  \end{eqnarray*}
  (arrow in $\Lambda X$).
\end{Exa}
For any $k\in \ZZ$, we have then the 2-functor
\begin{eqnarray*}
  \gpd & \longrightarrow & \gpd^{cen}\\
  X & \longmapsto & (\Lambda X,\xi^k).
\end{eqnarray*}
\begin{Exa} In the global quotient case, we have
  \begin{eqnarray*}
    \Lambda(M\mmod G)  \simeq  \coprod_{[g]}M^g\mmod C_g &\quad\text{  and
    } \quad&
    \xi^k\at{M^g} = g^k.
  \end{eqnarray*}
  Here $[g]$ is the conjugacy class, $M^g$ is the fixed point locus, and
  $C_g$ is the centralizer of $g$ in $G$. 
\end{Exa}
\begin{Def}\label{def:power_map}
  The $k$th power map 
  \begin{eqnarray*}
    \Pi^k\negmedspace: (\Lambda X,\xi^{k})&\longrightarrow&(\Lambda X,\xi^1)
  \end{eqnarray*}
  sends the object $g$ to $g^k$ and the arrow $h$ to $h$. We may also
  interpret $\Pi^k$ as a 1-morphism from $(\Lambda X,\xi^{km})$ to
  $(\Lambda X,\xi^{m})$. 
\end{Def}

\begin{Def}\label{def:kth-root}
  Let $X$ be a topological groupoid, let $\xi$ be an element of
  its center, and let $k$ be an integer. Then we define the groupoid
  \begin{eqnarray*}
    X[\xi^{\frac1k}] & := & \(\pt\mmod \ZZ\)\times \Lambda(X) /\sim,
  \end{eqnarray*}
  where the equivalence relation $\sim$ is generated by $k\sim\xi$. 
\end{Def}
\begin{Def}
  For an object $(X,\xi)$ of $\gpd^{cen}$, the subring 
  \begin{eqnarray*}
    K_{rot}(X,\xi)&\subset& K_{orb}(X)\ps{q^\frac1{|\xi|}}     
  \end{eqnarray*}
  is the Grothendieck group of formal power series
  \begin{eqnarray}\label{eq:F(q)}
    F(q) & = & \sum_{a\in\QQ_{\geq 0}} V_aq^q
  \end{eqnarray}
  satisfying the 
    {\em rotation condition} with respect to $\xi$:
    \begin{center}
      for each $a\in\QQ_{\geq 0}$ the coefficient
      $V_a$ is an $e^{2\pi ia}$-eigenbundle of $\xi$.
    \end{center}
  Here we are using Moerdijk's definition of orbifold $K$-theory
  $K_{orb}(X)$, see \cite[5.4]{Moerdijk}. 
\end{Def}
Power series satisfying the rotation condition may be thought of as
(infinite dimensional) vector bundles over 
\begin{eqnarray*}
  \RR\ltimes_\xi X & := & \(\pt\mmod \RR\)\times X\slash \sim 
\end{eqnarray*}
with $\sim$ generated by $1\sim\xi$.
\begin{Def}\label{def:K_Tate}
  Let $X$ be a topological groupoid. Then the Tate $K$-theory of $X$ is
  defined as
  \begin{eqnarray*}
    K_{Tate}(X) & := & K_{rot}(\Lambda(X),\xi^1).
  \end{eqnarray*}
\end{Def}
\subsection{Motivation}
\subsubsection{The Tate curve}
Let $S = \on{spec}(\ZZ\ps q)$.
Then there is an isomorphism of formal schemes over $S$
\begin{eqnarray*}
  \on{spf}\( K_{Tate}\(\CC P^\infty\) \) & = & {Tate(q)\widehat{{}_S}},
\end{eqnarray*}
making $K_{Tate}$ an elliptic cohomology theory, c.f.\
\cite[2.6]{Ando:Hopkins:Strickland:eswgtc}.
Let $C_k:=\ZZ/k\ZZ$ be the cyclic group with $k$ elements. Then
the complex representation ring of $C_k$ is 
\begin{eqnarray*}
  R(C_k)& \cong & \ZZ[\zeta_k],
\end{eqnarray*}
where the $k$th root $\zeta_k$ of one is identified with the
representation where $1$ acts as $e^{2\pi i/k}$.
So,
\begin{eqnarray*}
    spec\medspace\! R(C_k)\medspace &\cong& \mu_k,
\end{eqnarray*}
is the scheme of $k$th roots of unity over $\ZZ$.
where $\zeta_k$ is the irreducible representation with $\zeta_k(1) =
e^{2\pi i/k}$, and
one expects\footnote{See \cite{Devoto}, \cite{Ando:sigma}
  \cite{Greenlees}, \cite{Ginzburg:Kapranov:Vasserot}, \cite{Ganter:moonshine}}
\begin{eqnarray}
  \label{eq:level}
  \on{spec}\( K_{Tate}(\pt\mmod C_k) \) & = & {Tate(q)[k]},
\end{eqnarray}
where the right-hand side is the scheme of points of order $k$ in $Tate(q)$. 
By \cite[VII(1.13)]{Deligne:Rapoport}, 
\begin{eqnarray*}
  Tate(q)[k] & = & spec(D_k),
\end{eqnarray*}
where
$D_k$ is the Hopf-ring
\begin{eqnarray*}
  \label{eq:D_k}
  D_k & := & 
  \bigoplus_{j=0}^{k-1}\ZZ\ps q[x,x\inv]/x^k-q^j.
\end{eqnarray*}
After inverting $q$, we have the map
\begin{eqnarray*}
J\negmedspace:  \ZZ(\!(q)\!)[x] & \longrightarrow & \ZZ[\zeta_k](\!(q^\frac1k)\!)\\\
   q& \longmapsto & q\\
   x&\longmapsto  & \zeta_kq^{\frac jk}.
\end{eqnarray*}
The map $J$ 
identifies the $j$th summand of $D_k$ with the 
Laurent series satisfying the rotation condition for the group element $j$.
So, \eqref{eq:level} holds over the locus $spec \medspace\ZZ(\!(q)\!)$, where
the Tate curve is non-singular. 
The ring $D_k$ is closely related to isogenies into
$Tate(q)$ and plays a key role in Ando's definition of power operations
on $K_{Tate}$, see \cite[p.26 ff]{Ando:poweroperations}. 

We would also like to draw the reader's attention to Rezk's
\cite{Rezk:lectures}, where a 
much deeper connection between the constructions in \cite{Deligne:Rapoport} and
our Definition \ref{def:K_Tate} is discussed.
\subsubsection{Constant loops} 
The original motivation for Definition \ref{def:K_Tate} came from the
theory of orbifold loop spaces.
In \cite{Lupercio:Uribe}, Lupercio and
Uribe identify the inertia groupoid as the full subcategory 
\begin{eqnarray*}
  \Lambda X & \sub & \mF\!un (\RZ,X)
\end{eqnarray*}
of functors that are constant on objects.
Viewing the right-hand side as part of the loop groupoid
\begin{eqnarray*}
  LX & = & {\mO rb}\(\SS^1,X\),
\end{eqnarray*}
then leads them to the identification 
\begin{eqnarray*}
  \Lambda X & = & (LX)^\RR
\end{eqnarray*}
of the inertia groupoid with the orbifold loops that are fixed by
the rotation action. 
Under this identification, $\xi^k$ agrees with the rotation action of
$k\in\ZZ$ on (constant) orbifold loops.
\begin{Exa}[global quotients]
  By \cite[4.1]{Lupercio:Uribe}, 
  \begin{eqnarray*}
    L\(M\mmod G\) & \simeq & \coprod_{[g]}\mL_gM
  \end{eqnarray*}
  is the union of the {\em``twisted loop spaces''} of $M$, i.e., of paths
  in $M$ from $x$ to $gx$.
  This is what motivated Devoto's definition of
  equivariant Tate $K$-theory \cite[pp.5f]{Devoto}. 
\end{Exa}
In this context, the rotation condition is motivated by the following fact.
\begin{Lem}
  Let $V$ be a finite dimensional, real orbifold vector bundle on
  $X$. Then we have a dense isomorphism
  \begin{eqnarray*}
    LV\at{\Lambda(X)} & \cong & V_0 \oplus\bigoplus_{a\in\QQ_+}V_a^\CC
    q^a.
  \end{eqnarray*}
Here $V_0$ is the summand of $V\at{\Lambda(X)}$ fixed by $\xi^1$, and
for $a\in\QQ_+$, the bundle $V_a^\CC$ is the $e^{2\pi ia}$-eigenbundle
of $\xi^1$ inside the complexification $V_{\Lambda(X)}^\CC$, while
$q^a$ indicates that $\RR$ acts on this summand via $t\mapsto e^{2\pi iat}$.
\end{Lem}
\begin{proof}
  Let $(x,g)$ be an object in $\Lambda(X)$. By
  \cite[4.1.1]{Lupercio:Uribe}, we may identify $(LV)_{(x,g)}$ with the space
  of loops
  $$
    \{
      \map\gamma\RR V_x\mid \gamma(t+1) = g\gamma(t)
    \}.
  $$
  So, 
  \begin{eqnarray*}
    LV_{(x,g)} &  \sub & maps\(\RR/|g|\ZZ,V_x\).
  \end{eqnarray*}
  The Fourier expansion principle gives a dense isomorphism 
  \begin{eqnarray*}
    maps\(\RR/|g|\ZZ,V_x\)  & \cong & V_{x,0}\oplus\bigoplus_{n\geq 1}V_x^\CC
               q^{\frac n{|g|}},
  \end{eqnarray*} 
  under which
  $LV_{(x,g)}$ maps to the submodule specified in the
  the Lemma. 
\end{proof}
\subsubsection{Moonshine}\label{sec:moonshine}
Let $M$ be the monster, and consider the
central extension $\widetilde\Lambda_\alpha\(\pt\mmod M\)$ of the
inertia groupoid $\Lambda\(\pt\mmod M\)$ classified by the Moonshine
cocycle 
$$
  \alpha\in H^3\(BM, \on{U(1)}\)
$$
(see \cite[Meta Thm p.29]{Mason}).
Let $\widetilde\xi$ be a lift of $\xi^1$ to an element of the center
of $\widetilde\Lambda_\alpha\(\pt\mmod M\)$. Then
\begin{eqnarray*}
  \widetilde\Lambda_\alpha\(\pt\mmod M\) & \simeq & \pt\mmod\widetilde C_g,
\end{eqnarray*}
where the $\widetilde C_g$ are central extensions of the centralizers
in $M$. A choice of $\widetilde\xi$ amounts to a choice of lift
$\widetilde g\in\widetilde C_g$ for each $g$ in a system of
representatives for the conjugacy classes of $M$. 
A Laurent series
\begin{eqnarray*}
  F & \in & \bigoplus_{n\geq 0} \mR\!ep\(\widetilde
  C_g\)\!\(\!\!\({q^{\frac1{|\widetilde g|}}}\)\!\!\).
\end{eqnarray*}
satisfies the rotation condition
if and only if its character, a function $F\(g,h;q^{\frac1{|\widetilde g|}}\)$
in $h\in\widetilde C_g$, satisfies
\begin{eqnarray*}
  F\(g,\widetilde gh;q^{\frac1{|\widetilde g|}}\) & = & 
  F\(g,h; e^{\frac{2\pi i}{|\widetilde g|}}q^{\frac{1}{|\widetilde g|}}\).
\end{eqnarray*}
Interpreting the coefficients as characters of projective
representations of $C_g$, this is
\begin{eqnarray*}
  F(g,gh;\tau) & = & \zeta\cdot F(g,h;\tau+1), 
\end{eqnarray*}
where $\zeta$ is a root of unity and $q=e^{2\pi i \tau}$. These are
Condition (3) and half
of Condition (1) of Norton's generalized Moonshine conjecture \cite{Norton}. 
\subsection{Properties}\label{sec:properties}
Tate $K$-theory is an {\em orbifold theory with
  transfers}. In other words, it satisfies the following list of properties.
\subsubsection{Orbifold theory}
Let $$\longmap f{(X,\xi)}{(Y,\upsilon)}$$ be a 1-morphisms in
$\gpd^{cen}$. Then $f$ induces a ring map 
$$
  \longmap{f^*}{K_{rot}(Y,\upsilon)}{K_{rot}(X,\xi)}.
$$
If $f$ and $g$ are naturally isomorphic then $f^*=g^*$, and if $f$ is
an equivalence of groupoids then $f^*$ is an isomorphism. It follows
that the analogous statements hold for 
$K_{Tate}$ and that the latter is a well-defined functor on Moerdijk's
orbifold category.
\subsubsection{Sums and products }
We have natural isomorphism
%
%
\begin{eqnarray*}
  K_{Tate}(\emptyset) &\cong&   \{0\},\\
  K_{Tate}(\pt) & \cong & \ZZ\ps q,\\
  K_{Tate}(X\sqcup Y) &\cong&   K_{Tate}(X) \oplus  K_{Tate}(Y),
\end{eqnarray*}
and natural maps
\begin{eqnarray*}
  K_{Tate}(X)\tensor K_{Tate}(Y) & \longrightarrow &
  K_{Tate}(X\times Y), 
\end{eqnarray*}
making $K_{Tate}$ a lax rig functor from $Gpd$ to abelian rings.

\medskip
\noindent{\bf Exactness:}
Fix a finite group $G$. Then the functor
\begin{eqnarray*}
  M&\longmapsto & K_{Tate}(M\mmod G)
\end{eqnarray*}
sends cofiber sequences of $G$-spaces to exact sequences.
\begin{proof}
Indeed, for any element $g\in G$, the functor
\begin{eqnarray*}
  G\text{-spaces} & \longrightarrow &  C_g\text{-spaces} \\
  X & \longmapsto & X^g
\end{eqnarray*}
preserves cofiber sequences.
\end{proof}
\subsubsection{Transfers}
In this section all groupoids are assumed to have finite stabilizer groups. 
Given two maps of groupoids
$$
  X\xrightarrow{\phantom{xx}u\phantom{xx}} Z
  \xleftarrow{\phantom{xx}a\phantom{xx}} Y,
$$
we can form the {\em fiber} square
\begin{equation}
  \label{eq:fibre_square}
  \xymatrix{%
  X\times_ZY \ar[rr]^{v}\ar[dd]_b&& Y\ar[dd]^a\\ \\ 
  X\ar[rr]^u && Z,
  }
\end{equation}
as in \cite[2.3]{Moerdijk}, commuting up to a natural isomorphism $\eta$
(part of the data) and universal with respect to this
property. 
We will refer to the groupoid $X\times_ZY$ as the {\em
  fibred product} of $X$ and $Y$ over $Z$. Its objects are triples
$(x,y,u(x)\xrightarrow{g}a(y))$, consisting of an object of each $X$ and $Y$ and an arrow
between their images in $Z$. Arrows between two such triples $(x,y,g)$ and
$(x',y',g')$ are 
are pairs of arrows $(x\xrightarrow{h}x',y\xrightarrow{k}y')$ in $X$
and $Y$
satisfying $g'u(h) = a(k) g$.
Examples are the 
{\em comma category} groupoids
$y\uparrow f$ and $Y\uparrow f$, defined by the fiber squares
$$
  \xymatrix{%
  y\uparrow f \ar[rr]^{p_X}\ar[dd]&& X\ar[dd]^f\\ \\ 
  \pt\ar[rr]^y && Y,
  }
  \quad\quad\quad  \quad
  \xymatrix{%
  Y\uparrow f \ar[rr]^{P_X}_\sim\ar[dd]_{P_Y}&& X\ar[dd]^f\\ \\ 
  Y\ar@{=}[rr] && Y.
  }
$$ 
%
\begin{Def}
  We say that $\map fXY$ is {\em essentially a finite cover} if for each
  object $y\in Y$, there is a neighborhood $U$ of $y$ in $\ob Y$ 
  such that
  \begin{eqnarray*}
    U\times_YX & \simeq & U\times (y\uparrow f),
  \end{eqnarray*}
  and $y\uparrow f$ is equivalent to a finite groupoid.
\end{Def}
This condition ensures that pull-back of vector bundles possesses a right-adjoint
$$
  \longmap{RKan_f}{Vect_\CC(X)}{Vect_\CC(Y)},
$$
the {\em right Kan extension} along $f$, with
\begin{eqnarray*}
  \(RKan_fV\)_y & = & \lim_{\longleftarrow}{}_{y\uparrow f}\(p_X^*V\).
\end{eqnarray*}
%

\begin{Def}\label{def:!}
  Let $\map f{X}{Y}$ be essentially a finite
  cover. Then we write 
  $$
    \longmap{f_!}{K_{orb}(X)}{K_{orb}(Y)}
  $$
  for the map induced by $LKan_f$.
  We will refer to $f_!$ as {\em transfer along f}.
\end{Def}
\begin{Lem}
  Let $\longmap f{(X,\xi)}(Y,\upsilon)$ be a 1-morphism in
  $\gpd^{cen}$. If $\xi$ acts with eigenvalue $e^{2\pi i a}$ on the
  vector bundle $V$ on $X$ then $\upsilon$ acts with the same
  eigenvalue on $RKan_fV$. 
\end{Lem}
\begin{proof}
  The action of $\upsilon$ on
  \begin{eqnarray*}
    (RKan_fV)_y & \sub & \bigoplus_{(x,g)} V_x
  \end{eqnarray*}
  is by the permutation of the summands sending the object
  $(x,g)$ of $y\uparrow f$ to $(x,\upsilon_y g)$. The limit condition
  applied to the arrow $\map{\xi_x}{(x,\upsilon_yg)}{(x,g)}$ in 
  $y\uparrow f$ forces an element
  of the limit to satisfy
  $$
    v_{(x,g)} \quad = \quad \xi_x(v_{(x,\upsilon_y g)}) \quad=\quad
    e^{2\pi i a}\cdot v_{(x,\upsilon_y g)}.
  $$
  This proves the claim.
\end{proof}
\begin{Cor}
  We have transfers along essentially finite
  covers in $K_{rot}$ and in $K_{Tate}$.   
\end{Cor}
We will denote these also by $f_!$.   
So, $f_!$ in $K_{Tate}$ stands for $(\Lambda f)_!$ in $K_{rot}$.

\begin{Prop}\label{prop:push-pull}
  Given a fiber square as in \eqref{eq:fibre_square},
  assume that the map $a$ is essentially a finite cover. Then $b_!$ is
  defined, and we have
  \begin{eqnarray*}
    u^* a_! & = & b_! \medspace v^*.
  \end{eqnarray*}
\end{Prop}
\begin{proof}
  Let $x$ be an object of $X$. Then we have an equivalence of groupoids
  \begin{eqnarray*}
    x\uparrow b & = & \pt\times_X(X\times_Z Y)\\
                & \simeq & \pt\times_ZY \\
                & = & u(x)\uparrow a.
  \end{eqnarray*}
  Consider the canonical natural transformation of functors
  $\on{Vect_\CC}Y \to \on{Vect_\CC}X$
  \begin{eqnarray*}
    u^*RKan_a& \Longrightarrow & RKan_b v^*.
  \end{eqnarray*}
  Restricted to the fibers, this
  is the composition of isomorphisms
  \begin{eqnarray*}
    \(RKan_bv^* V\)_x & = & \lim_{\longleftarrow}{}_{x\uparrow
      b}\(p_{(X\times_ZY)}^*v^*V\)\\
      &\cong & \lim_{\longleftarrow}{}_{u(x)\uparrow a}\(p_Y^*V\)\\
      & = & \(RKan_a\)_{u(x)}.
  \end{eqnarray*}
  Further, $\Lambda$ preserves fiber squares.
\end{proof}
\begin{Cor}
  Orbifold Tate $K$-theory, restricted to the 2-category $\gpd^{fin}$
  of finite groupoids, is a global Mackey functor.
\end{Cor}
Here we are using the definition of global Mackey functor spelled out
in \cite{Ganter:2-lambda}. It is well known that, in the same manner,
orbifold $K$-theory gives the Mackey functor sending a finite
group(oid) to its representation ring.
\subsection{Characters}
\subsubsection{$n$-class functions}
Let $X$ be a finite groupoid, let $R$ be a ring, and let $n\geq 0$ be
a natural number.
\begin{Def}
  An $n$-class function on $X$ with values in $R$ is an $R$-valued map
  \begin{eqnarray*}
    \chi\negmedspace :[\Lambda^n X] & \longrightarrow & R,
  \end{eqnarray*}
  defined on the set of isomorphism classes of the $n$-fold inertia
  groupoid of $X$.
  We denote the ring of all such maps by
  \begin{eqnarray*}
    \on{n-Class}(X,R) & \cong & H^0(\on{Borel}(\Lambda^nX);R). 
  \end{eqnarray*}  
\end{Def}
Explicitly, a $0$-class function is a function on the isomorphism
classes of $X$, and for $n\geq 1$, an $n$-class function is
defined on $n$-tuples of commuting automorphisms of $X$ and satisfies
\begin{eqnarray*}
  \chi(g_1,\dots,g_n) & = &  \chi(sg_1s\inv,\dots,sg_ns\inv).
\end{eqnarray*}
Here all the $g_i$ are automorphisms of the same object $x$.
It is well-known that $\on{0-Class}(-,R)$ 
is a global Mackey functor, namely group cohomology,
with transfers along faithful maps $\map fXY$. It follows that
$\on{n-Class}(-,R)$ is also a Mackey functor, whose transfers along
faithful maps are given by
\begin{eqnarray}
  \label{eq:0transfer}
  f_!(\chi)(g_1,\dots,g_n)& = & \sum_{[\ul h]\mapsto[\ul
    g]}\frac{|aut(\ul g)|}{|aut(\ul h)|}\cdot \chi(h_1,\dots,h_n).
\end{eqnarray}
Here $\ul g$ abbreviates $(g_1,\dots, g_n)$, and $[\ul g]$ and
$aut(\ul g)$ are, respectively, its
isomorphism class and its automorphism group in $\Lambda^nY$, and similarly for $\ul
h$.

If $R$ is a $\QQ$-algebra, then \eqref{eq:0transfer} makes sense for
all $f$ and extends $\on{n-Class}(-;R)$ to a Mackey
functor with all transfers. 
\begin{Exa}
  Let $\map jHG$ be an inclusion of groups and write $j$ also for the
  corresponding map of groupoids $\map j{\pt\mmod H}\pt\mmod G$. Then
  \begin{eqnarray*}
    j_!(\chi)(g_1,\dots,g_n) & = & \frac1{|H|}\sum_{s\ul g s\inv\in
      H^n}\chi(sg_1s\inv,\dots,sg_ns\inv). 
  \end{eqnarray*}
\end{Exa}
\begin{Exa}\label{exa:eps_G}
  Let $G$ be a finite group, and let $\map {\varepsilon_G}{\pt\mmod G}\pt$ be the
  unique map. Then
  \begin{eqnarray*}
    (\varepsilon_G)_!(\chi)(1) & = & \frac1{|G|}\sum_{\ul g}\chi(g_1,\dots,g_n), 
  \end{eqnarray*}
  where the sum runs over all $n$-tuples of commuting elements of $G$.
\end{Exa}
\subsubsection{Character theory}
Write $\gpd^f$ for the 2-category of finite groupoids, and 
assume that $E$ is an orbifold theory with transfers as in Section \ref{sec:properties}.
\begin{Def}
  We say that $E$ has a
  {\em character theory} 
  if there exists a ring $R$ and a natural transformation of Mackey functors
  \begin{eqnarray*}
    \chi\negmedspace:E\at{\mG pd^{f}} & \Longrightarrow & \on{n-Class}(-,R)
  \end{eqnarray*}
  for some $n$. 
    We will refer to $n$ as the {\em height} of the
  theory.
\end{Def}
Here $R$ is assumed to be either a $\QQ$-algebra or torsion free with
the understanding that some transfers in
$\on{n-Class}(-,R)$ are, a priori, only defined over $\QQ\tensor R$.
Naturality of $\chi$ then implies 
integrality of transfers in the image of $\chi$.
\begin{Exa}
  $K$-theory has a character theory of height $1$ with
  \begin{eqnarray*}
    R & = & \on{\underrightarrow{lim}}_k\medspace \ZZ[\zeta_k],
  \end{eqnarray*}
  where $\zeta_k$ is a $k$th root of $1$.  
  Note that the
  universal level $k$-structure on the multiplicative group has its
  home over $\on{spec}(\ZZ[\zeta_k])$.
  The integrality statement above amounts to the well-known fact that 
  \begin{eqnarray*}
    \frac1{|G|}\sum_{g\in G}\chi(g) & \in & \ZZ
  \end{eqnarray*}
  when $\chi$ is the character of a representation. 
\end{Exa}
\begin{Exa}
  Tate $K$-theory has a height 2 character theory, taking 
  values in  
  \begin{eqnarray*}
    \on{\underrightarrow{lim}}_k\medspace\(\ZZ[\zeta_k]\ps{q^\frac1k}\).
  \end{eqnarray*}
  By \cite[VII.2.4,VII.1.16.4]{Deligne:Rapoport}), 
  \begin{eqnarray*}
    \on{spec}\ZZ[\zeta_k]\ps{q^\frac1k}& \longrightarrow & \(\mM_k\)
  \end{eqnarray*}
  is the completion at infinity of the level $k$ stack $\mM_k$ defined
  in \cite[p.36, IV.3.2,3.5]{Deligne:Rapoport}.
\end{Exa}
\begin{Exa}
  Fix a prime $p$, and
  let $E_n$ be Borel equivariant Morava-Lubin-Tate theory.
  This has a height $n$ character theory (\cite[Thms.\ C\&D]{Hopkins:Kuhn:Ravenel},
  \cite[7.9]{Ganter:thesis}), taking values in 
  \begin{eqnarray*}
    L(E_n) & = & S^{-1} \on{\underrightarrow{lim}}_k\medspace E_n(B(\ZZ/p^k\ZZ)^n)
  \end{eqnarray*}
  where the multiplicative set $S$ consists of the Euler classes of
  non-trivial representations. The ring 
  \begin{eqnarray*}
    D_k & := & im\(E_n(B(\ZZ/p^k\ZZ)^n)\to S\inv E_n(B(\ZZ/p^k\ZZ)^n)\)
  \end{eqnarray*}
  is the home of the universal level $k$-structure on
  the formal group of $E_n$ (see \cite[Thm 3.3.2]{Ando:isogenies}).
\end{Exa}
The character theory in \cite{Hopkins:Kuhn:Ravenel} is a lot more
powerful than our summary here suggests. It would be interesting to
formulate Hopkins-Kuhn-Ravenel theory for orbifolds and to 
see exactly how much of the Hopkins-Kuhn-Ravenel story can be told for
$K_{Tate}$. It seems reasonable to expect that $K_{Tate}$ fits into
Stapleton's framework of transchromatic character maps \cite{Stapleton}.
\subsubsection{A formula for Induction}
Let $H\sub G$ be an inclusion of finite groups, $M$ a $G$-manifold,
and let $$\longmap j{M\mmod H}{M\mmod G}$$ be the inclusion. Then the
transfer $j_!$ in $K_{orb}$ is
$$
  \longmap{ind_H^{\phantom{.}G}}{K_H(M)}{K_G(M)}.
$$  
We will write 
$$
  \longmap{I_H^{\phantom{}G}}{K_{Tate}(M\mmod H)}{K_{Tate}(M\mmod G)}
$$
for the transfer $j_!$ in $K_{Tate}$.
\begin{Prop}
  Let $a$ be an element of $K_{Tate}(M\mmod H)$. Then $I_H^G(a)$ is the
  element of $$\bigoplus_{[g]}K_{C_G(g)}(M^g)\ps{q^{\frac1{|g|}}}$$
  whose $[g]$th summand equals
  \begin{eqnarray*}
    \label{eq:I}
    I_H^G(a)_{[g]} & = &
    \begin{cases}
      ind_{C_H}^{\phantom{.}C_G}\(a_{[g]}\) & \text{ if $g\in H$}\\
      0 & \text{ if $[g]\cap H = \emptyset$.}
    \end{cases}
  \end{eqnarray*}
\end{Prop}
\begin{proof}
  Because of the commuting diagram
  $$
    \xymatrix{%
    M\mmod H\ar[0,2]^j && M\mmod G \\
    (G\times_HM)\mmod G, \ar[-1,0]^\sim\ar[-1,2]_{pr_2}
    }
  $$
  $I_H^G$ may be identified with
  the $G$-equivariant (Atiyah) transfer in $K_{rot}$ along 
  $\Lambda(pr_2)$. This may be viewed as the map of $G$-spaces
  \begin{eqnarray*}
    G\times_H ob(\Lambda(M\mmod H)) & \longrightarrow &
    ob(\Lambda(M\mmod G)\\
    (s,(x,h)) & \longmapsto & (sx,shs\inv).
  \end{eqnarray*}
  Hence we have
  \begin{eqnarray*}
    I_H^G(a)_{[g]} & = & \sum_{r\inv gr}r\cdot a_{[r\inv gr]},
  \end{eqnarray*}
  where $r$ runs over a system of representatives of $G/H$. The set of
  representatives contributing to the sum is given by the image of 
  the inclusion 
  $$
    C_G/C_H \hookrightarrow G/H,
  $$
  so we may as well choose the $r$s as a set of representatives of
  $C_G/C_H$. 
\end{proof}
%
%
\section{Power operations}
\subsection{Symmetric powers of orbifolds}
Recall that $Gpd$ 
is a symmetric bimonoidal category (``rig'' category) with the
monoidal structures given by $(\sqcup,\emptyset)$ and
$(\times,\pt)$. Every groupoid $X\in ob(\mG pd)$ is a
monoid with respect to $(\sqcup,\emptyset)$, via the fold map
$$\longmap d{X\sqcup X}X.$$ 
A groupoid $X$ is a monoid with
respect to $(\times,\pt)$ if and only if $X$ is endowed with the 
structure of a symmetric monoidal category.
\begin{Def}
  Then $n$th symmetric power of $X\in ob(Gpd)$ is the groupoid
  $S_n\smallint X$ obtained from $X^n$ by adding the additional arrows
  $$
    (x_1,\dots,x_n) \xrightarrow{\phantom{xx}\sigma\phantom{xx}}
    (x_{\sigma\inv(1)},\dots,x_{\sigma\inv(n)})\quad\quad\sigma\in S_n,
  $$
  composing with the arrows in $X^n$ as follows:
  \begin{eqnarray*}
    \sigma\circ (g_1,\dots,g_n) & = &
    (g_{\sigma\inv(1)},\dots,g_{\sigma\inv(n)})\circ\sigma .
  \end{eqnarray*}
  The total symmetric power of $X$ is
  \begin{eqnarray*}
    S(X) := \coprod_{n\geq 0}S_n\smallint X.
  \end{eqnarray*}
\end{Def}
The endofunctor $S$ of $Gpd$ is exponential in the sense that it is
a monoidal functor
$$
  \longmap S{(Gpd,\sqcup,\emptyset)}{(Gpd,\times,\pt)}.
$$
It follows that $S(X)$ comes equipped with a symmetric monoidal
structure 
\begin{eqnarray*}
  * := S(d),
\end{eqnarray*}
which turns out to be concatenation.
The unit of $*$ is the unique object  $()$ of $S_0\smallint X$. The
triple 
$$
  (S(X),*,())
$$
may be viewed as the free symmetric monoidal category on $X$. More precisely,
the functor
\begin{eqnarray*}
  X&\longmapsto& (S(X),*,())
\end{eqnarray*}
is the left adjoint to the forgetful functor from the category of
monoids in $Gpd$ to $Gpd$.
In particular, $S$ is a monad.
\begin{Def}\label{def:monad}
  We will write
  $$
    \mu\negmedspace : S^2 \Longrightarrow S
  $$
  and 
  $$
    \iota\negmedspace : \id \Longrightarrow S
  $$
  for the structure maps of the monad $S$. Explicitly, these are given
  by the inclusions $S_n\smallint S_m \sub S_{nm}$ and by
  the inclusion of $X=
  S_1\smallint X$ inside $S(X)$.
\end{Def}
\subsection{The inertia groupoid of $S(X)$}
Objects of $\Lambda(S(X))$ are given by quadruples $(n,\ul x,\ul
g,\sigma)$, where $n\geq 0$ is a natural number, $\ul x\in ob(X^n)$,
and $\sigma\in S_n$, while $\ul g$ is a system of arrows
$x_i\xrightarrow{g_i}x_{\sigma(i)}$. Note that this notation is
somewhat redundant, all the information is contained in the pair
$(\sigma;\ul g)$. The inertia groupoid inherits a monoidal structure
from $S(X)$, given by
\begin{eqnarray*}
  (\sigma;g_1,\dots,g_n) * (\tau;h_1,\dots,h_m) & = &
  (\sigma\sqcup\tau;g_1,\dots,g_n,h_1,\dots,h_m).
\end{eqnarray*}
An object $(\sigma,\ul g)$ of $\Lambda(S(X))$ is indecomposable with
respect to $*$ if and only if $n>0$ and $\sigma$ acts on
$\{1,\dots,n\}$ with a single orbit (i.e., $\sigma$ is a long cycle).
Otherwise, the cycle decomposition of $\sigma$ yields a decomposition
of $(\sigma,\ul g)$ into indecomposables.
\begin{Def}
\label{def:Phi} 
  Let $\varsigma_k\in S_k$ be the cycle $\varsigma_k=(1,\dots,k)$, and let
  $\Phi_k(X)$ be the full subgroupoid of $\Lambda(S(X))$ with objects
  $(\varsigma_k,\ul x)$.
  Let $\Phi(X)$ be the groupoid
  \begin{eqnarray*}
    \Phi(X) &:=& \coprod_{k\geq 1}\Phi_k(X).
  \end{eqnarray*}
  Further, let $\varphi\in\cen(\Phi_k)$ be the restriction of
  $\xi^1$ to $\Phi(X)$.
\end{Def}
The essential image of $\Phi(X)$ inside $\Lambda(SX)$ is the subgroupoid
of indecomposable objects.
The functor $\Phi$ is additive. i.e.,
$\Phi(\emptyset) = \emptyset$, and $\Phi$ preserves $\sqcup$. 
\begin{Lem}\label{lem:Phi}
  We have an equivalence of monoidal groupoids
  $$
    \longmap Q{S(\Phi(X))}\Lambda(S(X)),
  $$
  which is natural in $X$ and satisfies 
  \begin{eqnarray*}
    Q S(\varphi) & = & \xi^1Q.
  \end{eqnarray*}
\end{Lem}
\begin{proof}
  Let $I$ be the inclusion
  $$
    I\negmedspace : \Phi(X) \hookrightarrow \Lambda(S(X)),
  $$
  and let $Q$ be the composite
  $$
    Q\negmedspace : S(\Phi(X))
    \xrightarrow{\phantom{xx}S(I)\phantom{xx}}
    S(\Lambda(S(X))) \xrightarrow{\phantom{xx}\varepsilon\phantom{xx}}\Lambda(S(X)),
  $$
  where the second map is the counit of the adjunction $(S,*,())\dashv
  forget$. 
  Then $Q$ is monoidal. Since the essential image of $I$ consists
  exactly of the indecomposable objects of $\Lambda(S(X))$, it follows
  that $Q$ is essentially surjective. One checks that $Q$ is also fully
  faithful. (This boils down to the fact that $\sigma$ and $\tau$ are
  conjugate in $S_n$ if and only if their cycle decompositions are
  congruent.) 
\end{proof}

\begin{Lem}\label{lem:E_k}
  For each $k\geq 1$, there is an equivalence of groupoids
  \begin{eqnarray*}
    E_k\negmedspace : \Phi_k(X)     & \xrightarrow{\phantom{xx}\sim\phantom{xx}} &
    \Lambda(X)[\xi^{\frac1k}], 
  \end{eqnarray*}
  identifying $\varphi$ with $\xi^{\frac1k}$.
\end{Lem}
\begin{proof}
  Let 
  $(\ul x,\ul g)$ be an object of $\Phi_k(X)$. For $1\leq i\leq k$, 
  we let $\widehat g_i\in\ob(\Lambda(X))$ be the
  composite 
  $$
    \widehat g_i\negmedspace:\quad
    x_i\xrightarrow{\phantom{i}g_i\phantom{i}}
    x_{i+1}\xrightarrow{\phantom{i}g_{i+1}\phantom{i}} 
    x_{i+2}\xrightarrow{\phantom{i}g_{i+2}\phantom{i}}\dots
    x_{i-1}\xrightarrow{\phantom{i}g_{i-1}\phantom{i}}x_i.
  $$
  Then $E_k$
  sends $(\ul x,\ul g)$ to the object $\widehat g_1$ of
  $\Lambda(X)$. So, $E_k$ is surjective on objects.
  We have
  \begin{eqnarray}
  \label{eq:phi^k}
  \varphi_{(\ul x,\ul g)}^k & = & (id, \ul {\widehat g}).
  \end{eqnarray}
  Let $(\varsigma_k^m,\ul h)$ be 
  an arrow 
  in $\Phi_k(X)$. 
  Then $(\varsigma_k^m,\ul h)$ can be factored as 
  \begin{eqnarray*}
     (id,\ul h')\circ \varphi^m& = &
     \varphi^m \circ  (id,\ul h')
  \end{eqnarray*}
  with $m\in \ZZ$, and this presentation is unique up to the relation
  \eqref{eq:phi^k}. 
  Now the map $E_k$ sends the arrow
  $(id,\ul h)$ to $h_1$ and the center element $\varphi$ to $\xi^{\frac1k}$.
  
  Consider an arrow of the form $(id,\ul h)$ from $(\ul x,\ul g)$ to
  $(\ul x',\ul g')$  
  in $\Phi_k(X)$. Then the $h_i$ fit into a commuting diagram
  \begin{equation}
    \label{eq:Phi_arrow}
    \xymatrix{
    x_1 \ar[0,2]^{g_1}\ar[d]_{h_1}& &  x_2 \ar[0,2]^{g_2}\ar[d]_{h_2} &&  x_3
    \ar[r]\ar[d]_{h_3} & \cdots\ar[r]&
    x_k \ar[0,2]^{g_k}\ar[d]_{h_k}& &  x_1 \ar[d]^{h_1} \\
    x'_{1} \ar[0,2]^{g'_{1}} &&  x'_{2} \ar[0,2]^{g'_{2}} &&  x'_{3}
    \ar[r] & \cdots\ar[r]&
    x'_{k} \ar[0,2]^{g'_{k}} &&  x'_{1}. 
    }
  \end{equation}
  In particular, 
  \begin{eqnarray}
    \label{eq:conj}
    h_1 \widehat g_1 h_1\inv & = & \widehat g_{1}',
  \end{eqnarray}
(and similarly for $i$). Hence $E_k$ is well-defined. Further,
the $h_i$ are uniquely determined by $h_1$, and
any $h_1$ satisfying \eqref{eq:conj} gives rise a to system $\ul h$
fitting into \eqref{eq:Phi_arrow}. Hence $E_k$ is fully faithful.
\end{proof}
\begin{Def}
  Let $F_k$ be the quasi-inverse of $E_k$ that sends the object $g$ to 
  $$(\varsigma_k;g,id,\dots,id)$$
  ($k-1$ times $id$).
\end{Def}
All the above constructions are Morita invariant.
\subsubsection{Comparison to Dijkgraaf Moore Verlinde Verlinde} 
Let $X$ be the groupoid of a global quotient orbifold $M\mmod G$,
where $G$ is a finite group. Then the loop space of $S_n\smallint X$,
\begin{eqnarray*}
  L(S_n\smallint X) & \simeq & \mF\!un(\RZ,S_n\smallint X),
\end{eqnarray*}
has objects $(\sigma;\ul g,\ul\gamma)$ where $(\sigma;\ul g)\in S_n\smallint
G$ and $\ul \gamma$ is an $n$-tuple of paths
$$
  \longmap{\gamma_i}\RR M
$$
satisfying
\begin{eqnarray}
\label{eq:path}
  g_i\gamma_i(t) & = & \gamma_{\sigma(i)}(t+1).
\end{eqnarray}
\begin{Exa}
  Assume that $\sigma = \varsigma_k$. Then we have
  \begin{eqnarray*}
    \gamma_1(t+k) & = & g_k\cdots g_1\gamma_1(t)\\
    \gamma_2(t+k) & = & g_1g_k\cdots g_2\gamma_2(t)
  \end{eqnarray*}
  and so on, and each $\gamma_i$ determines the others via
  \eqref{eq:path}.
  Hence $$\gamma_i\in\mF\!un(\RR\mmod k\ZZ, X)$$ may be thought of as
  a loop of length $k$ in $M\mmod G$.
\end{Exa}
Let $L_kX\sub L(S_n\smallint X)$ be the full subgroupoid whose objects
have $\sigma=\varsigma_k$. Then our groupoid $\Phi_k$ 
can be identified with $L_kX^\RR$, the part of $L_kX$
that is invariant under the rotation action by $\RR$. So, $\Phi_kX$
may be viewed as the groupoid of constant loops of length $k$ in
$M\mmod G$, and $\varphi$ is identified with the rotation action by
$1\in\ZZ$ on these long loops.
The key argument in \cite{Dijkgraaf:Moore:Verlinde:Verlinde} is
summarized by the equivalence of groupoids
\begin{eqnarray}\label{eq:DMVV}
  LS(X) &\simeq& 
  S\(\coprod_{k\geq 1}L_kX\)
\end{eqnarray}
(see also \cite{Tamanoi}). 
Our Lemma \ref{lem:Phi} above follows from \eqref{eq:DMVV} by
restricting to constant loops on both sides. 
\subsection{Power operations in orbifold theories}
Let $E$ be an orbifold theory with products. 
\begin{Def}
  A {\em total power operation} for $E$ is   
  a (non-linear) natural transformation
  $$
    P\negmedspace : E \Longrightarrow E\circ S
  $$
  satisfying
  \begin{description}
    \item[Comodule property:] $P$ makes $E$ a comodule over the comonad
      $(-)\circ S$. 
    \item[Exponentiality:] The map $$\longmap P{E(\emptyset)}E(\pt)$$ sends
      $0$ to $1$, and 
    $$
      \longmap{P}{E(X\sqcup Y)}{E(SX\times SY)}
    $$
    sends $(a,b)$ to the external product $P(a)\tensor P(b)$. 
  \end{description}
  We can write $P$ as 
  $P=(P_n)_{n\geq 0}$, with $P_n:=\iota_n^* P$. Then $P_n$ is called the
  $n$th power operation of $P$. 
\end{Def}
\begin{Exa}[Atiyah Power Operations]
\label{exa:Atiyah}
  Let $V$ be an orbifold vector bundle on $X$. Then we have the Atiyah Power
  Operations
  \begin{eqnarray*}
    P_n(V) & := & S_n\ltimes V^{\tensor n}
  \end{eqnarray*}
  of \cite{Atiyah}.
  We will see below how to extend this definition to virtual vector
  bundles. For now we note that $P=(P_n)_{n\geq 0}$ satisfies the
  axioms of a total power operation, whenever it is defined. 
\end{Exa}
\subsubsection{Consequences of the definition}
\begin{Prop}
  \begin{enumerate}
  \item We have $P_1 (a) = a$. 
  \item More generally, 
    \begin{eqnarray*}
      i_k^*P_k(a) = a^{\tensor k} 
    \end{eqnarray*}
    (external product), where $i_k$ is the inclusion of $X^k$ in
    $S_k\smallint X$. 
  \item We have
    \begin{eqnarray*}
      (PS)\circ P = \mu^*\circ P, 
    \end{eqnarray*}
    or, equivalently,
    \begin{eqnarray*}
      P_n(P_m(a)) = res^{S_{nm}}_{S_n\smallint S_m}\( P_{nm}(a)\).
    \end{eqnarray*}
  \item $P_0(a) = 1$.
  \item The external product of $a\in E(X)$ and $b\in E(Y)$ is given
    by
    \begin{eqnarray*}
      a\tensor b  &=& i^* P(a,b), 
    \end{eqnarray*}
    where $(a,b)\in E(X\sqcup Y)$, and $$i:=(\iota_X\times\iota_Y)
    :X\times Y \hookrightarrow SX\times SY$$ is product of the
    canonical inclusions (see Definition \ref{def:monad}).
  \item We have
    \begin{eqnarray*}
      P(a\tensor b) & = & j^*\(P(a)\tensor P(b)\)
    \end{eqnarray*}
    (external products), where 
    $$
      j\negmedspace :S(X\times Y)\longrightarrow
      S(X)\times S(Y),
    $$
    is defined as
    \begin{eqnarray*}
      j=(S(\pi_X),S(\pi_Y)),
    \end{eqnarray*}
    with $\pi_X$ and $\pi_Y$ the projections to the respective factors. 
  \item The map $P$, and equivalently all the $P_n$, preserves interior products:  
    \begin{eqnarray*}
      P_n(ab) & = & P_n(a)P_n(b),
    \end{eqnarray*}
    for $a,b\in E(X)$.
  \item For each $n$, we have $P_n(1)=1$.
  \end{enumerate}
\end{Prop}
\begin{proof}
  The Comodule Property translates into 
  Points 1 and 3. To prove Point 2, we apply naturality of $P$ to the
  fold map $\map{f}{X\sqcup \dots \sqcup X}X$, noting that $S(f)$ is
  the concatenation product $*$. This leads to a commuting diagram
  $$%
    \xymatrix{
      E(X)\ar[0,2]^P\ar[d]_{f^*}&& E(SX) \ar[d]^{(*)^*}\ar[0,2]^{(\iota_ki_k)^*}&& E(X^k)\\
      E(X\sqcup\dots\sqcup X)\ar[0,2]^P&& E((SX)^k),
      \ar[-1,2]_{(\iota_1^k)^*} 
    }
  $$
  where $f^*(a) = (a,\dots,a)$. By the Exponential Property, the bottom
  composite sends this to $P_1(a)^{\tensor k}$, hence 2. 

  For 4, we apply naturality of $P$ to the unique map $\map e\emptyset
  X$, arriving at
  \begin{eqnarray*}
    \iota_0^*\circ P\at X & = & P\at\emptyset\circ e^*.
  \end{eqnarray*}
  The left-hand side of this is $P_0$. By the first part of the
  Exponential Property, the right-hand side is the
  constant map $1$.
  
  5 is an immediate consequence of the Exponential Property and 1.

  To prove 6, one checks that $j=\mu\circ S(i)$ and considers the diagram
  $$
    \xymatrix{%
      E(X\sqcup Y) \ar[0,2]^P\ar[d]_P && E(SX\times SY) \ar[0,2]^{i^*}\ar[d]^P &&
      E(X\times Y)\ar[d]^P\\
      E(SX\times SY) \ar[0,2]_{\mu^*} && E(S^2(X\sqcup Y)) \ar[0,2]_{(Si)^*}& &
      E(S(X\times Y)),
    }
  $$
  whose left square commutes by 3, and whose right square commutes by
  naturality of $P$. By 5, the composite on the top row sends $(a,b)$
  to $a\tensor b$. By the Exponential Property, the leftmost vertical
  arrow sends $(a,b)$ to $P(a)\tensor P(b)$.
  
  7 follows from 6 by considering the special case $X=Y$: Let
  $\map\delta XX\times X$ be the diagonal map.
  \begin{eqnarray*}
    P(ab) & = & P\(\delta^*\(a\tensor b\)\) \\
          & = & (S\delta)^* P(a\tensor b) \\
          & = & (S\delta)^*j^* (P(a)\tensor P(b))\\
          & = & \delta_{S(X)}^* \(P(a)\tensor P(b)\)\\
          & = & P(a)P(b).
  \end{eqnarray*}
  Finally, we show 8. In the special case $X=\pt$, this follows from
  3. and the first part of the Exponential Property. The general case
  follows by applying naturality of $P$ to the unique map from $X$ to $\pt$.
\end{proof}
\subsection{The graded ring $E(S(X))$}
From now on, we assume that $E$ possesses transfers along (faithful)
essentially finite covers. 
Then transfer along the concatenation product
$$
  \longmap{*}{S(X)\times S(X)}{S(X)}
$$
defines a second multiplication $\bullet$ on the ring $E(S(X))$,
making it a graded ring with respect to the grading
\begin{eqnarray*}
  E(S(X)) & \cong & {\bigoplus_{n\geq 0}^\wedge}\medspace E(S_n\smallint X)
\end{eqnarray*}
with unit ${\bf 1}$ in degree zero. 
\begin{Prop}
  Assume that $P$ is natural also with respect to transfers, and let
  $a,b\in E(X)$.
  Then we have  
  \begin{eqnarray*}
    P(0) & = & {\bf 1}\quad\text{and}\\
    P(a+b) & = & P(a) \bullet P(b).
  \end{eqnarray*}
\end{Prop}
\begin{proof}
  The first equality follows by applying naturality of $P$ to the
  transfer along the inclusion of
  $\emptyset$ in $X$. The second follows by applying naturality of $P$
  to the transfer along the fold map $X\sqcup X\to X$. 
\end{proof}
In the global quotient case, all the relevant transfers exists by
equivariant Spanier-Whitehead duality, and $P$ is automatically
natural with respect to transfers obtained in this manner.
It seems reasonable to expect that, in a suitably defined category of
`orbispectra', a similar argument would make the demand for naturality
along faithful transfers redundant and that all $K$-theoretic
transfers above should be induced by maps in a
$K_{orb}$-localization of this orbispectra category. 
\begin{Cor}
  The Atiyah power operations of Example \ref{exa:Atiyah} can be extended, in a
  unique way, to virtual vector bundles.
\end{Cor}
\begin{proof}
  Indeed, $P$ takes values in ${\bf 1}+ K_{orb}(SX)$ and hence $P([V])$ is
  invertible with respect to $\bullet$ for any vector bundle $V$.
\end{proof}
Similarly, one argues that the definition of Atiyah power operations
extends to $K_{rot}$, giving 
$$
  \longmap P{K_{rot}(X,\xi)}{K_{rot}(SX,S\xi)}.
$$
\subsection{Power operations in $K_{Tate}$}
Recall from Lemma \ref{lem:Phi} that we have
\begin{eqnarray}
\label{eq:KS}
  K_{Tate}(S(X)) & \cong & K_{rot}\(S(\Phi(X)), S(\varphi)\).
\end{eqnarray}
This becomes a graded isomorphism if, 
on the right-hand side, elements supported on $S_n\smallint\Phi_k(X)$ are
given degree $nk$.
\begin{Def}\label{def:s_k}
  Let $X$ be an orbifold groupoid, let $\xi$ be an element of its
  center, and let $k\geq 1$ be a natural number.
  Then we define the map
  \begin{eqnarray*}
    s_k\negmedspace : K_{rot}(X) & \longrightarrow &
    K_{rot}(X[\xi^{\frac1k}])\\
    \left[\sum V_aq^a\right]&\longmapsto & \left[ \sum V_aq^{\frac ak} \right],   
  \end{eqnarray*}
  where $\xi^{\frac1k}$ acts on the coefficient $V_a$ by $e^{2\pi ia/k}$.
\end{Def}
  The map $s_k$ commutes with the Atiyah power operations in the
  following sense:
  let $(X,\xi)$ be as in Definition \ref{def:s_k}, recall the notation
  $X^{(k)}$ for $X[\xi^{\frac1k}]$ and let 
  \begin{eqnarray*}\label{eq:s_kP}
    \delta\negmedspace : (SX)^{(k)} & \hookrightarrow & 
    S(X^{(k)}) 
  \end{eqnarray*}
  be the canonical
  inclusion (i.e., $\delta$ sends $(S\xi)^{\frac 1k}$ to $S(\xi^{\frac 1k})$).
  Then we have
  \begin{eqnarray}
    s_k \circ P & = & \delta^* P s_k.
  \end{eqnarray}
\begin{Def}
  Let 
  $$
    \longmap{\theta}{K_{Tate}(X)}{K_{rot}(\Phi(X),\varphi)}
  $$
  be the additive operation whose $k$th component is $E_k^*\circ
  s_k$. Here $E_k$ is the equivalence defined in Lemma \ref{lem:E_k}.
\end{Def}
\begin{Def}
  The total power operation in $K_{Tate}$ is defined as the composite
  $$
    P^{Tate}\negmedspace : K_{Tate}(X)
    \xrightarrow{\phantom{xx}\theta\phantom{xx}}
    K_{rot}(\Phi X,\varphi) \xrightarrow{\phantom{xx}P\phantom{xx}}
    K_{rot}(S\Phi X,S\varphi) \xrightarrow{\phantom{xx}(Q^*)\inv\phantom{xx}}
    K_{Tate}(SX). 
  $$
\end{Def}
\begin{Thm}
  This $P^{Tate}$ satisfies the axioms of a total power operation.
\end{Thm}
\begin{proof}
  The exponential property follows immediately from that of $P$ and
  from additivity of $\theta$,
  \begin{eqnarray*}
    \theta\negmedspace : K_{Tate}(X\sqcup Y) &
    \longrightarrow &
    K_{rot}(\Phi X\sqcup \Phi Y)\\
    (a,b) & \longmapsto & (\theta a,\theta b).
  \end{eqnarray*}

  \begin{Lem}
    Let the natural transformation $\vartheta$ be defined as the composite
    $$
      \vartheta_X\negmedspace: \Phi SX \xrightarrow{\phantom{xx}I_{SX}\phantom{xx}}
      \Lambda S^2X\xrightarrow{\phantom{xx}\Lambda\mu_X\phantom{xx}} \Lambda SX.
    $$
    Then we have 
    \begin{eqnarray*}
      \(\theta S\)\circ P^{Tate} & = & \vartheta^* P^{Tate}.
    \end{eqnarray*}
  \end{Lem}
  \begin{proof}
    It suffices to prove this after restricting both sides 
    of the equation to $\Phi_kS$. 
    Consider the commuting diagram
    $$%
      \xymatrix{
      \Lambda^{(k)}S  \ar[0,2]^{F_kS}_\simeq && \Phi_kS\ar[0,2]^\vartheta&&
      \Lambda S\\
      \\
      \(S\Phi\)^{(k)}\ar[-2,0]^{Q^{(k)}}_\simeq\ar[0,2]^\delta &&
      S(\Phi^{(k)})\ar[0,2]^{Su}_\simeq& &S\Phi\ar[-2,0]^Q_\simeq
      }
    $$
    where the equivalence $u$ sends the object $(\varsigma_l;\ul g)$
    of $\Phi_l^{(k)})$ to
    the object
    $$
      (\varsigma_{kl};\underbrace{id,\dots,id}_{k-1},g_1,\underbrace{id,\dots,id}_{k-1},g_2,id,
      \dots,\dots,id,g_l),  
    $$
    of $\Phi_{kl}$,
    sends a morphism of the form $(id,\ul h)$ to
    $$
      (id;\underbrace{h_1,\dots,h_1}_{k},\underbrace{h_2,\dots,h_2}_{k},\dots,
      \dots,\underbrace{h_l,\dots,h_l}_{k}),   
    $$
    and $\varphi^{\frac1k}$ to $\varphi$.
    We have
    \begin{eqnarray*}
      u^*\theta & = & s_k \theta
    \end{eqnarray*}
    and hence 
    \begin{eqnarray*}
       (Q^{(k)})^*\circ F_k^* \circ\vartheta^*\circ P^{Tate} 
      & = & \circ\delta^*\circ(Su)^*\circ P\circ\theta\\
      & = & \circ\delta^*\circ P\circ u^*\circ \theta\\
      & = & \circ\delta^*\circ P\circ s_k\circ\theta\\
      & = & s_k \circ P\circ \theta\\
      & = & (Q^{(k)})^*\circ s_k\circ  P^{Tate}\\
    \end{eqnarray*}
    Since $(F_k\circ Q^{(k)})^*$ is an isomorphism the claim of the
    Lemma follows. 
  \end{proof}
  To complete the proof of the theorem, we consider the commuting diagram
  $$%
    \xymatrix{
      S\Phi S\ar[0,2]^{S\vartheta}\ar[2,0]_{QS}^\simeq && S\Lambda
      S\ar[2,0]^\varepsilon &&
      S^2\Phi\ar[0,-2]_{SQ}^\simeq\ar[2,0]^{\mu\Phi}\\ \\
      \Lambda S^2 \ar[0,2]^{\Lambda\mu} && \Lambda
      S && S\Phi\ar[0,-2]_{Q}^\simeq,
    }
  $$
  where the middle map, $\varepsilon$, is the counit of the adjunction
  $(S,*,())\dashv forget$ (compare to the definition of $Q$ in the
  proof of Lemma \ref{lem:Phi}). 
  This yields
  \begin{eqnarray*}
    (P^{Tate}S)\circ P^{Tate} 
    & = & ((QS)^*)\inv\circ (P\Phi S)\circ
    \vartheta^*\circ(Q^*)\inv\circ(P\Phi)\circ \theta\\ 
    & = & ((QS)^*)\inv\circ
    (S\vartheta)^*\circ((SQ)^*)\inv\circ(PS\Phi)\circ(P\Phi)\circ
    \theta\\ 
    & = & ((QS)^*)\inv\circ
    (S\vartheta)^*\circ((SQ)^*)\inv\circ (\mu\Phi)^*\circ(P\Phi)
    \circ \theta\\ 
    & = & (\Lambda\mu)^*(Q^*)\inv\circ(P\Phi) \circ \theta\\ 
    & = & (\Lambda\mu)^*P^{Tate}. 
  \end{eqnarray*}
\end{proof}
\subsection{Symmetric and exterior powers}
Let $E$ be a theory with transfers and power operations. Write
$$
  \longmap{\varepsilon}{S(\pt)}\mathbb N
$$
for the augmentation map, sending $S_n\smallint\pt$ to $n$.
We may identify $E(\mathbb N\times X)$ with the formal power series
ring $E(X)\ps t$. 
\begin{Def}
  The total $E$-theoretic symmetric power $S_t$ is the composite
  $$
    S_t\negmedspace : EX  \xrightarrow{\phantom{xx}P\phantom{xx}}
    ESX   \xrightarrow{\phantom{xx}\delta^*\phantom{xx}}
    E(S(\pt)\times X)
    \xrightarrow{\phantom{xx}\epsilon_!\phantom{xx}} E(X)\ps t,
  $$
  where $\delta$ is the diagonal map, sending the arrow $(\sigma,g)$
  to the arrow $(\sigma;g,\dots,g)$ with $n$ copies of $g$ if $\sigma\in S_n$.
\end{Def}

\begin{Prop}\label{prop:ring-maps}
  In the definition of $S_t$, the maps $\delta^*$ and $\varepsilon_!$
  are maps of graded rings.
\end{Prop}
\begin{Cor}
  The total symmetric power is exponential:
  \begin{eqnarray*}
    S_t(a+b) & = & S_t(a)\cdot S_t(b).
  \end{eqnarray*}
\end{Cor}
\begin{proof} We now prove the proposition.
The ring multiplication on $E(X)\ps t$ corresponds to
$$
  \xymatrix{
    E(\mathbb N\times X\times\mathbb N\times X)
    \ar[0,2]^{\delta_X^*}&&
    E(\mathbb N\times X\times\mathbb N) \ar[0,2]^{m_!}&&
    E(\mathbb N\times X), 
  }
$$
where $\map{\delta_X}XX\times X$ is the diagonal map and $\map
m{\mathbb N\times\mathbb N}{\mathbb N}$ is multiplication.
Similarly, one defines the graded multiplication on $E(S(\pt)\times
X)$ as the push-pull 
$$
  E\(\(S(\pt)\times X\)^2\)
  \xrightarrow{\phantom{xx}{\delta_X^*}\phantom{xx}}
  E\(S(\pt)^2\times X\) 
  \xrightarrow{\phantom{xx}{*_!}\phantom{xx}}
    E(S(\pt)\times X).
$$
The claim then follows by a tedious but straight-forward, iterated
application of Proposition \ref{prop:push-pull}.
\end{proof}
\subsection{Symmetric powers in Tate $K$-theory}
\begin{Def}
  For a positive natural number $k$, we define the operator
  \begin{eqnarray*}
    \beta_k\negmedspace : K_{Tate}(X) & \longrightarrow &
    K_{Tate}(X)\\
    F&\longmapsto&\left[\Pi_k^*F(q^\frac1k)\right]^{\ZZ/k\ZZ},
  \end{eqnarray*}
  where $\Pi_k$ is the $k$th power map of Definition
  \ref{def:power_map}, and the generator of $\ZZ/k\ZZ$ acts on the
  coefficient $\Pi_k^*(V_a)$ of $q^{\frac ak}$ by $\xi^1\cdot
  e^{-2\pi i\frac ak}$. So, the $\ZZ/k\ZZ$-invariant part is
  the largest possible subspaces of the coefficients such
  that the result satisfies the rotation condition for $\xi^1$.   
\end{Def}
The action of the $\beta_k$ should be compared to the action of
$\widehat\ZZ^\times$ in \cite{Hopkins:Kuhn:Ravenel}.

Say the order of $\xi^1$ is $r$ and we are given an element
$F(q^\frac1r)\in K_{Tate}(X)$. A priori, it looks like $\beta_k(F)$ is
a power series in $q^\frac1{kr}$, but the rotation condition implies
that it is in fact a power series in $q^\frac1r$.  
\begin{Exa}
  If $X=M$ is a manifold then $F$ is of the form $F=\sum_{n=0}^\infty
  V_nq^n$, and 
  \begin{eqnarray*}
    \beta_k(F) & = & \sum_{m=1}^\infty V_{km}q^m.
  \end{eqnarray*}
  This should be compared to
  \cite[(2.13)]{Dijkgraaf:Moore:Verlinde:Verlinde}. 
\end{Exa}
\begin{Thm}\label{thm:total-symmetric-power}
  The total symmetric power in Tate $K$-theory is described by the
  formula
  \begin{eqnarray}\label{eq:Witten-class}
    S_t^{Tate}(F) & = & \bigotimes_{k=1}^\infty S_{t^k}(\beta_k(F)).
  \end{eqnarray}
  Here $S_t$ stands for the total Atiyah symmetric power on $K_{Tate}$.
\end{Thm}
\begin{Cor}
  If $X=M$ is a manifold and $F=V$ is a vector bundle on $M$, then
  $S_t^{Tate}(V)$ is Witten's stable exponential characteristic class
  of $V$.
\end{Cor}

\begin{proof}[Proof of the Theorem]
We will construct a commuting diagram:

\begin{center}
\includegraphics{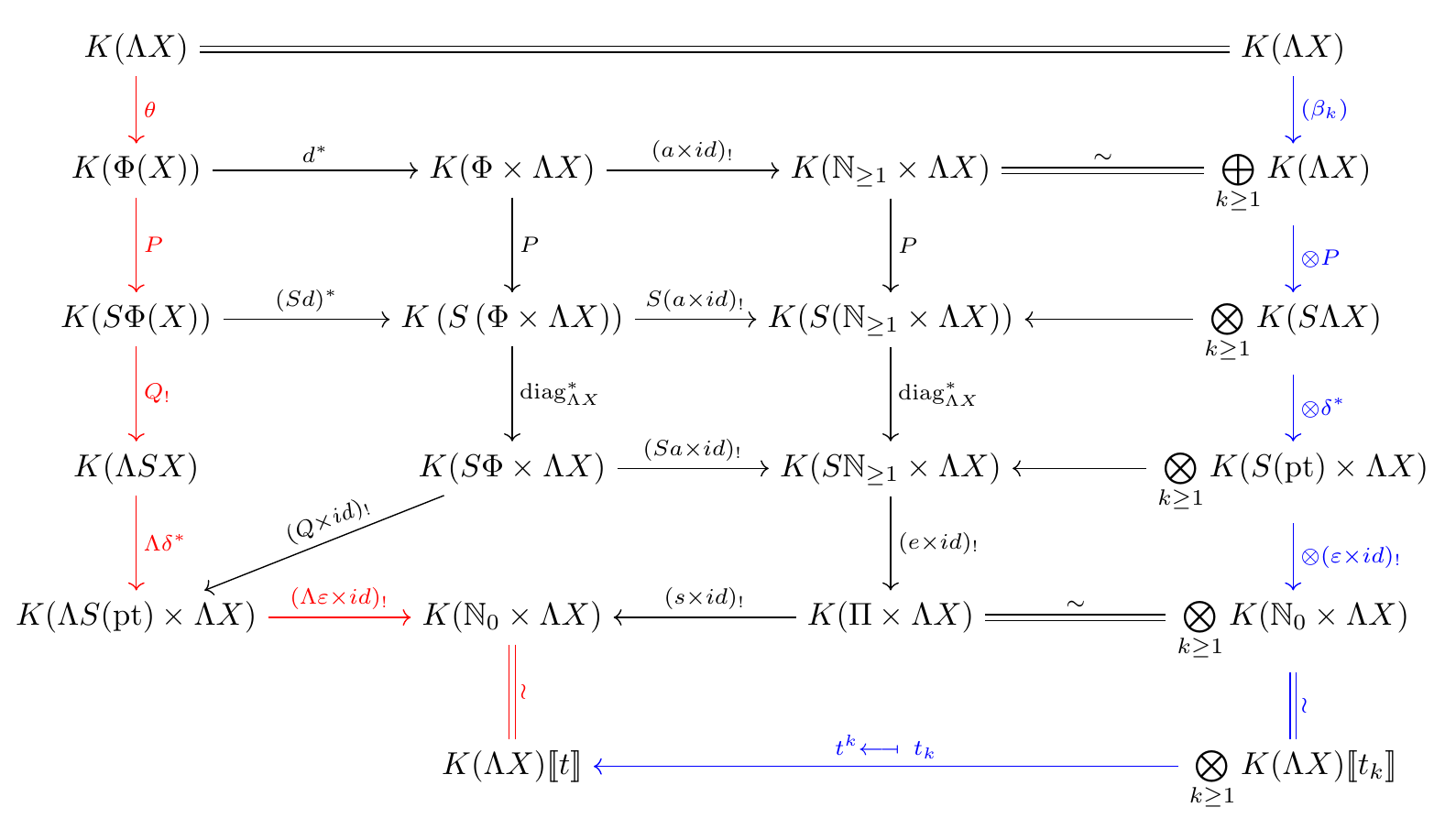}
\end{center}

Again, $K$ is short for $K_{rot}$.
The composite of the red arrows gives the left-hand side of
\eqref{eq:Witten-class}, that of the blue arrows gives the right-hand
side.
The other arrows are defined as follows: abbreviating $\Phi(\pt)$ as
$\Phi$, we let 
$$\map d{\Phi\times\Lambda X}{\Phi(X)}$$
be the restriction of $\Lambda\delta$ to the source of $d$. 
So, the $k$th component $d_k$ of $d$
sends the object $g$ to $(\varsigma_k;g,\dots,g)$ and the
morphism
$(\varsigma_k^m,h)$ to $(\varsigma_k^m;h,\dots,h)$.
Let $$\map{a}{\Phi}\NNo$$ be the map that sends $\Phi_k(\pt)$ to
$\{k\}$.
Together with $\theta$, these maps form the top square of the diagram.

The three commuting squares in the second row are almost immediate
from the properties of $P$: the right-most follows, by
induction over the degree, from the exponential property.
The unlabeled horizontal maps out of the tensor products on the very
right all are obtained as 
the (internal) product of projection maps of the form 
\begin{eqnarray*}
  {S\(\coprod_k Y_k\)}&\longrightarrow &{S(Y_j)}, 
\end{eqnarray*}
mapping $(y_i)$ to the list of entries from $Y_j$.

Next, we have the maps $\on{diag_{\Lambda X}}$ defined as
\begin{eqnarray*}
{\on{diag_{\Lambda X}}}\negmedspace: {S\NNo\times\Lambda
  X}&\longrightarrow&{S(\NNo\times\Lambda X)} \\ 
((k_1,\dots,k_n),g) & \longmapsto &
((k_1,g),\dots,(k_n,g)) \quad\text{on objects}\\
(\sigma_n,h) &
\longmapsto &
(\sigma_n;h,\dots,h)) \quad\quad\quad\text{on morphisms,}
\end{eqnarray*}
and similarly, with $\Phi$ in the role of $\NNo$.

Finally, we write $\Pi$ for the set of isomorphism classes of $\Lambda
S(\pt)$ or, equivalently, $S\Phi$ or $S\NNo$. Elements of $\Pi$ may be
thought of as sequences $(n_k)_{k\geq 1}$ with all but finitely many
entries equal to zero, and we have the degree map
\begin{eqnarray*}
  s\negmedspace : \Pi& \longrightarrow & \NN_0\\
  (n_k)_k & \longmapsto& \sum_kkn_k, 
\end{eqnarray*}
and the canonical quotient map
\begin{eqnarray*}
  e\negmedspace : S\NNo&\longrightarrow &\Pi,
\end{eqnarray*}
fitting into the pentagon in the fourth row.
The right square in the same row is obtained from the family of
commuting squares

\begin{center}
  \includegraphics{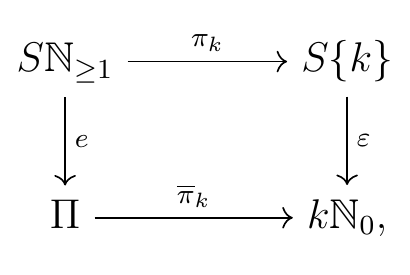}
\end{center}
where $\ol\pi_k$ is the map of isomorphism classes induced
by $\pi_k$. The fact that $e_!\circ\prod\pi_k^*$ equals $\bigotimes
\epsilon_!$ is not immediate. The argument boils down to the fact
that, for a $G$-representation $V$ and an $H$-representation $W$, the
invariant part $(V\tensor W)^{G\times H}$ equals $V^G\tensor W^H$.
\end{proof}
\subsection{Adams operations and Hecke operators}
The total symmetric power in $K$-theory has the well-know generating
function
\begin{eqnarray}\label{eq:Adams-operators}
  S_t(x) & = & \exp\(\sum_{m\geq 1}\frac{\psi_m(x)}mt^m\),
\end{eqnarray}
where the $\psi_m$ are the {\em Adams operations}. In fact,
\eqref{eq:Adams-operators} is often used as a definition of the Adams
operators. 
As a consequence of Theorem \ref{thm:total-symmetric-power}, we obtain
the following generating function for $S^{Tate}_t$:
\begin{eqnarray}\label{eq:generating-function}
  S^{Tate}_t(F) & = & \exp\(\sum_{m\geq 1}T_m(F)t^m\),
\end{eqnarray}
where the {\em Hecke operators} $T_m$ are defined as
\begin{eqnarray*}
  T_m(F) & := & \frac1a\sum_{ad = m}\psi_a\(\beta_d(F)\).
\end{eqnarray*}
\begin{Exa}
  It is well known that the Adams operations on the coefficient rings
  $K_G(\pt)=R(G)$ have the following effect on characters:
  \begin{eqnarray*}
    \psi_m(\chi)(g) & = & \chi(g^m).
  \end{eqnarray*}
  Let $F$ be an element of the coefficient ring 
  \begin{eqnarray*}
    K_G^{Tate}(\pt) & \cong & \bigoplus_{[g]} R(C_g)\ps{q^\frac{1}{|g|}}.
  \end{eqnarray*}
  As in Section \ref{sec:moonshine}, we view $F$ as the $q$-expansion
  of a function 
  $F(g,h;\tau)$, with $q^\frac1{|g|} = e^{2\pi i\tau/|g|}$. Then
  \begin{eqnarray*}
    \beta_d(F)(g,h;\tau) & = & \frac1d\sum_{0\leq
      b<d}F\(g^d,g^{-b}h;\frac{\tau+b}d\),
  \end{eqnarray*}
  and hence
  \begin{eqnarray*}
    T_m(F)(g,h;\tau) & = & \frac1m\sum_{ad=m}\sum_{0\leq
      b<d}F\(g^d,g^{-b}h^a;\frac{\tau+b}d\).
  \end{eqnarray*}
  These are the equivariant Hecke operators 
  that play an important role in Moonshine, see 
  \cite{McKay}, \cite{Ganter:moonshine}, \cite{Morava} and \cite{Carnahan}. 
\end{Exa}
The equivariant Hecke operators can be defined by an equivariant Hecke
correspondence, and 
this is the sense in which our power operations on $K_{Tate}$ are {\em
  `elliptic'}. A natural question is whether our notion of elliptic
can be strengthened to resemble the definition of {\em $H_\infty$
  elliptic spectrum} in \cite[Def.16.4]{Ando:Hopkins:Strickland:sigma}. In
other words, do our power operations define {\em descent data for
  level structures} on the Tate curve as in
\cite{Ando:Hopkins:Strickland:sigma}, and do our power operations 
specialize to the ones Ando defines in \cite[6.3]{Ando:poweroperations}?
\subsection{Exterior powers and replicability}
The {\em total exterior power} $\Lambda_t^{Tate}$ is defined by the equality
\begin{eqnarray*}
  \Lambda^{Tate}_t(F) & = & \(S_{-t}^{Tate}(F)\)\inv,
\end{eqnarray*}
so
\begin{eqnarray*}
  \Lambda^{Tate}_{-t}(F) & = & \exp\(-\sum_{m\geq 1}T_m(F)t^m\).
\end{eqnarray*}
Let $$F(q)\in R(G)(\!(q)\!)$$ be of
the form 
\begin{eqnarray*}
  F(q)&=&q\inv+a_1 q+a_2 q^2+\dots.
\end{eqnarray*}
From the Moonshine literature,
such $F$ are known as McKay-Thompson series. We recall\footnote{Compare
  e.g.\ \cite[(2.1)]{Teo} with $b=1$, $t=1/z$ and $F(q) = g(z)$.}
that the Faber polynomials $\Phi_{m,F}$ of $F$ are defined by
$$
  -\sum_{m=1}^\infty \Phi_{m,F} (w) t^m  = \log\(t(F(t)-w)\).
$$
Hence $\Phi_{m,F}$ is a polynomial in $w$ of degree $m$, which depends
on the first $m$ coefficients of $F$ and is uniquely characterized by
the fact that 
it is of the form 
$$\Phi_{m,F}(F(q))=q^{-m}+b_1q+b_2q^2+\dots.$$
\begin{Def}
  Let $F$ be a McKay-Thompson series. We write $F^{(a)}$ for the $a^{th}$
  Adams operations applied to $F$. We call $F$ {\em replicable}, if for
  every natural number $m$, we have
  \begin{eqnarray*}    
    \Phi_{m,F}(F(q)) & = & \sum_{ad = m}\sum_{0\leq b<d}
    F^{(a)}\(\frac{a\tau +b}{d}\)\\
    & = & m\cdot T_m(F)(q)
  \end{eqnarray*}
  Here $q=e^{2\pi i\tau}$.
\end{Def}
This appears to be the right notion of replicability of McKay-Thompson
series, it is the one that turns up in \cite{Borcherds}.
It follows that a McKay-Thompson series $F$ is
replicable if and only if it satisfies 
\begin{eqnarray*}
  F(t)-F(q) &=& t\inv\cdot\Lambda^{Tate}_{-t}(F(q)).
\end{eqnarray*}
\bibliographystyle{alpha}
\bibliography{stringy.bib}
\end{document}